\documentclass[11pt]{amsart}
\usepackage{amssymb, color}
\usepackage{amsmath}
\usepackage{amsthm}
 \font\tenmsb=msbm10 at 12pt \font\sevenmsb=msbm7 at 8pt \font\fivemsb=msbm5 at
6pt
\newfam\msbfam
\textfont\msbfam=\tenmsb \scriptfont\msbfam=\sevenmsb \scriptscriptfont\msbfam=\fivemsb
\def\Bbb#1{{\tenmsb\fam\msbfam#1}}

\def\C{\Bbb C}

\def\N{\Bbb N}

\def\S{\Bbb S}

\begin{document}
\def \theequation{\arabic{section}.\arabic{equation}}
\newcommand{\reset}{\setcounter{equation}{0}}

\newcommand{\beq}{\begin{equation}}
\newcommand{\noi}{\noindent}
\newcommand{\eeq}{\end{equation}}
\newcommand{\dis}{\displaystyle}
\newcommand{\mint}{-\!\!\!\!\!\!\int}

\def \bx{\hspace{2.5mm}\rule{2.5mm}{2.5mm}} \def \vs{\vspace*{0.2cm}} \def
\hs{\hspace*{0.6cm}}
\def \ds{\displaystyle}
\def \p{\partial}
\def \O{\Omega}
\def \b{\beta}
\def \m{\mu}
\def \l{\lambda}
\def \le{\lambda^*}
\def \ul{u_\lambda}
\def \D{\Delta}
\def \d{\delta}
\def \s{\sigma}
\def \e{\varepsilon}
\def \a{\alpha}
\def \g{\gamma}
\def \R{\mathbb{R}}
\def \S{\mathbb{S}}
\def\qed{%
\mbox{ }%
\nolinebreak%
\hfill%
\rule{2mm} {2mm}%
\medbreak%
\par%
}
\newtheorem{thm}{Theorem}[section]
\newtheorem{lem}[thm]{Lemma}
\newtheorem{cor}[thm]{Corollary}
\newtheorem{prop}[thm]{Proposition}
\theoremstyle{definition}
\newtheorem{defn}{Definition}
\newtheorem{rem}[thm]{Remark}
\newenvironment{thmskip}{\begin{thm}\hfill}{\end{thm}}
\def \pr {\noindent {\it Proof:} }
\def \rmk {\noindent {\it Remark} }
\def \esp {\hspace{4mm}}
\def \dsp {\hspace{2mm}}
\def \ssp {\hspace{1mm}}

\title{Conformal metrics in $\R^{2m}$ with constant $Q$-curvature and arbitrary volume}
\author{Xia Huang}
\address{Department of Mathematics and Center for Partial Differential Equations, East China Normal
University, Shanghai 200062, P.R. China} \email
{xhuang1209@gmail.com}
\author{Dong Ye}
\address{IECL, UMR 7502, D\'epartement de Math\'ematiques, Universit\'e de Lorraine, B\^at. A, \^{\i}le de Saulcy, 57045 Metz, France} \email{dong.ye@univ-lorraine.fr}
 \subjclass{35J30, 53A30, 35J91} \keywords{conformal geometry, constant $Q$-curvature}

\date{}
\begin{abstract}
We study the polyharmonic problem $\Delta^m u = \pm e^u$ in $\R^{2m}$, with $m \geq 2$. In particular, we prove that {\sl for any} $V > 0$, there exist radial solutions of $\Delta^m u = -e^u$ such that
$$\int_{\R^{2m}} e^u dx = V.$$
It implies that for $m$ odd, given arbitrary volume $V > 0$, there exist conformal metrics $g$ on $\R^{2m}$ with positive constant $Q$-curvature and vol$(g) =V$.  This answers some open questions in Martinazzi's work \cite{M}.
\end{abstract}

\maketitle
\section{Introduction}
In $\R^{2m}$, if the conformal metric $g_u = e^{2u}|dx|^2$ satisfies $(-\Delta)^m u = Q(x)e^{2mu}$, it is well known that (see for instance \cite{C1}) $Q(x)$ is just the $Q$-curvature of the metric $g_u$, here $|dx|^2$ is the Euclidean metric.

\medskip
One interesting question in conformal geometry is to understand the metrics with constant $Q$-curvature, i.e.~to understand
solutions of $(-\Delta)^m u = Qe^{2mu}$ in $\R^{2m}$ with $Q \in \R$. Considering $v = u - \ln\lambda$ with $\lambda > 0$,
$(-\Delta)^m u = Qe^{2mu}$ is equivalent to $(-\Delta)^m v = \lambda^{2m}Qe^{2mv}$. Therefore the precise value of $Q$ is not important and we can reduce the study to $Q \in \{0, \pm 1\}$.

\medskip
 For $Q > 0$ case, let $(\S^{2m}, g_{\S^{2m}})$ be the standard unit sphere in $\R^{2m+1}$ and $\pi: \S^{2m} \rightarrow \R^{2m}$ be the stereographic projection, we know that $Q_{g_{\S^{2m}}} = (2m-1)!$ and
$$(\pi^{-1})^*g_{\S^{2m}} = \frac{4|dx|^2}{(1+|x|^2)^2} := g_\pi.$$
Consequently $Q(g_\pi) = (2m-1)!$. By scaling and translation, for any $x_0 \in \R^{2m}$, $\l > 0$,
\begin{align}
\label{01}
u_{x_0, \l} = \ln\frac{2\l}{1+\l|x - x_0|^2}
\end{align}
satisfies
\begin{align*}
Q\big(g_{u_{x_0, \l}}\big) = (2m-1)! \quad \mbox{and} \quad {\rm vol}\big(g_{u_{x_0, \l}}\big) = \int_{\R^{2m}} e^{2mu_{x_0, \l}} dx = {\rm vol}(\S^{2m}).
\end{align*}

For $m = 1$, Chen \& Li \cite{CL} proved that any solution to
\begin{align}
\label{02}
-\Delta u = e^{2u} \;\;\mbox{in }\R^2, \quad {\rm vol}(g_u) = \int_{\R^2} e^{2u} dx < \infty
\end{align}
is given by the formula \eqref{01}, hence ${\rm vol}(g_u) = {\rm vol}(\S^2)$, i.e.~any conformal metric $g$ in $\R^2$ with
positive constant Gauss curvature and finite volume is provided by the stereographic projection of $\S^2$ into $\R^2$.

\medskip
Remark that without the assumption of finite volume, Liouville \cite{Lio} showed that there are many other entire solutions to $-\Delta u = e^{2u}$ in $\R^2$.

\medskip
The situation is very different for $m > 1$. Consider the problem
\begin{align}
\label{03}
(-\Delta)^m u = (2m-1)! e^{2mu} \;\;\mbox{in }\R^{2m}, \quad {\rm vol}\big(g_u) = \int_{\R^{2m}} e^{2mu} dx < \infty,
\end{align}
Chang \& Chen \cite{CC} proved the existence of non-spherical solutions: When $m > 1$, for any $0 < V < {\rm vol}(\S^{2m})$,
there exists a solution to \eqref{03} such that vol$(g_u) = V$.

\medskip
The condition vol$(g_u) < {\rm vol}(\S^{2m})$ was not only suggested by technical reasons, but it is also necessary when $m =2$.
Indeed, let $m = 2$, Lin showed in \cite{L} that any solution to \eqref{03} verifies vol$(g_u) \leq {\rm vol}(\S^4)$, and
the equality holds if and only if the solution is spherical (i.e.~given by \eqref{01}). Moreover, when $m = 2$, Wei \& Ye
\cite{WY} proved the existence of solution with any asymptotic behavior at infinity predicted by Lin. In particular, it
means that for $m = 2$ and any volume $V$ less than ${\rm vol}(\S^{4})$, there exists a very rich family of non radial solutions of \eqref{03} with vol$(g_u) = V$. Thus the situation in $\R^4$ is somehow well understood for $Q > 0$.

\medskip
Recently, Martinazzi \cite{M} found striking and new phenomena for $m \geq 3$: The solutions to \eqref{03} can have volume
larger than ${\rm vol}(\S^{2m})$. More precisely,
\begin{itemize}
\item[(i)] for $m = 3$, there exists $V^* > 0$ such that for any $V \geq V^*$, we have a solution $u$ of \eqref{03} in $\R^6$ such that ${\rm vol}(g_u) = V$;
\item[(ii)] for any $m \geq 3$ odd, there exists $V_m > {\rm vol}(\S^{2m})$ such that for every $V \in (0, V_m]$, there is a solution $u$ of \eqref{03} satisfying ${\rm vol}(g_u) = V$.
\end{itemize}
However, he could not rule out that $V_3 < V^*$ in (i)-(ii) (when $m=3$) and he asked if a gap phenomenon is possible,
that is, could it be a volume $V$ in $(V_3, V^*)$ such that the problem \eqref{03} has no solution verifying ${\rm vol}(g_u) = V$?
He asked also if the result in (i) could be generalized for $m \geq 5$ odd.

\medskip
In this work, we generalize completely (i)-(ii) by proving that for $m \geq 3$ odd,
there exist solutions to \eqref{03} with arbitrary volume.
\begin{thm}
\label{odd}
For every $m\geq 3$ odd, and every $V \in(0, +\infty)$ , there admits a conformal metric in $\mathbb{R}^{2m}$ satisfying $Q_g\equiv (2m-1)!$ and ${\rm vol}(g) = V$.
\end{thm}

The result for $m$ even is less complete, but which still suggests that no gap phenomenon exists for \eqref{03}.
\begin{thm}
\label{even}
For every $m \geq 2$ even, let
$${\mathcal V} = \left\{ \int_{\R^{2m}} e^{2mu} dx, \;\mbox{ with radial function $u$ satisfying \eqref{03}} \right\}.$$
Then ${\mathcal V}$ is an interval.
\end{thm}

Our approach is to study respectively entire radial solutions of the following polyharmonic equations ($m \geq 2$):
\begin{equation}\label{1.1}
\Delta ^m u= -e^u\quad\text{in}~\R^{2m}
\end{equation}
and
\begin{equation}\label{1.1new}
\Delta ^m u= e^u\quad\text{in}~\R^{2m}.
\end{equation}
The main results are
\begin{thm}
\label{thm1.1}
Let $m \geq 2$. Then for any $V \in(0, +\infty)$, there exists a radial solution $u$ to \eqref{1.1} such that
$$\int_{\R^{2m}} e^u dx = V.$$
\end{thm}

\begin{thm}
\label{thm1.1new}
Let $m \geq 2$. If there exists an entire radial solution $u_0$ to \eqref{1.1new}, then for any $0 < V < \|e^{u_0}\|_{L^1(\R^{2m})}$, there exists a radial solution $u$ of \eqref{1.1new} such that
$$\int_{\R^{2m}} e^u dx = V.$$
\end{thm}

Notice that given a solution $u$ to \eqref{1.1} or \eqref{1.1new}, the function $v:= \frac{1}{2m}\left[ u - \ln(2m)!\right]$ solves
$$
(-\Delta)^m v = \pm(-1)^{m+1}(2m-1)!e^{2mv}\;\;\text{in}~\R^{2m},\quad
\int_{\R^{2m}} e^{2mv} dx = \frac{1}{(2m)!} \int_{\R^{2m}} e^u dx.
$$
Hence, Theorems \ref{odd} and \ref{even} are just direct consequence of Theorems \ref{thm1.1} and \ref{thm1.1new} respectively. So we
need just to prove Theorems \ref{thm1.1} and \ref{thm1.1new}.

\medskip
Furthermore, for the negative constant $Q$-curvature case, i.e.~when $Q < 0$, there is no entire solution to $\Delta u = e^{2u}$ in $\R^N$ for $m=1$ and any dimension $N \geq 1$ (see for example Theorem 1 in \cite{O}). Here again, we find a completely different situation for $m > 1$. Recently, Hyder \& Martinazzi showed
that for any $m \geq 2$, $V > 0$, and any polynomial $P(x)$ of degree $\leq (2m-2)$ verifying $\lim_{\|x\|\to\infty} x\!\cdot\!\nabla P(x) =  \infty$, there exists $u$ such that $(-\Delta)^m u = -(2m-1)!e^{2mu}$ in
$\R^{2m}$ and
$$\int_{\R^{2m}} e^{2mu} dx = V, \quad u(x) = -P(x) + \frac{2V}{{\rm vol}(\S^{2m})}\ln\|x\| + C + o(1) \quad \mbox{as }\|x\|\to\infty.$$
The above result is a direct consequence of Theorem 1.2 in \cite{HM} combined with Theorem C there, which was previously proved in \cite{M2}.

\section{Proof of Theorem \ref{thm1.1}}
\reset
\subsection{Preliminaries} Consider the following initial value problem in $\R^N$ for general dimensions $N \geq 3$.
\begin{equation}\label{2.4}
\begin{cases}
\begin{aligned}
&\Delta^m u=-e^u, && \\&u^{(2i+1)}(0)=0, &&  \forall\; 0\leq i\leq m-1,\\
&\Delta^iu(0)=a_i,&& \forall\; 0\leq i\leq m-1.
\end{aligned}
\end{cases}
\end{equation}
Here $u(x) = u(r)$ is a radial function, the laplacian $\Delta$ is seen as $\Delta u = r^{1-N}\left(r^{N-1}u'\right)'$
and $a_i$ are constants in $\R$. We will denote $u_{(a_i)}$ the radial solution to \eqref{2.4}.

\medskip
Clearly, there exist suitable constants $\alpha_i$ such that $\Phi_\alpha(r) = \sum_{0\leq j \leq m-1} \alpha_jr^{2j}$ verifies
\begin{align*}
\Delta^i\Phi_\alpha(0)=a_i, \quad \forall\; 0\leq i\leq m-1.
\end{align*}
As $\Delta^m(u_{(a_i)}-\Phi_\alpha) = -e^{u_{(a_i)}} < 0$, it's easy to check that $u_{(a_i)}(r) \leq \Phi_\alpha(r)$ whenever
$u_{(a_i)}$ exists. Therefore $-e^{u_{(a_i)}}$ is locally bounded whenever $u_{(a_i)}$ exists. Applying standard ODE theory, we can
claim that for any $(a_i) \in \R^m$, the unique radial solution of \eqref{2.4} is defined globally in $\R_+$, in other words,
an entire radial solution to $\Delta^m u = -e^u$ exists in $\R^N$ for any $(a_i)$.

\begin{rem}
\label{rem2.1}
For $N \ne 2m$, if $u$ is a solution to \eqref{2.4} with $e^u \in L^1(\R^N)$, we can get solution with arbitrary $L^1$ norm by the scaling $u_\l(x) = u(\l x) + 2m\ln\l$, since
\begin{align*}
\Delta^m u_\lambda = -e^{u_\lambda}, \quad \int_{\R^N} e^{u_\lambda} dx = \lambda^{2m-N}\int_{\R^N} e^u dx.
\end{align*}
So our main concern here is only relevant for $N = 2m$. We should mention that Farina \& Ferrero provide recently in \cite{FF}
many precise studies for radial solutions of $\Delta^m u = \pm e^u$ in $\R^N$ with general $m, N \in \N^*$.
\end{rem}

The following Lemma is inspired by \cite{FF}. It's a simple but important fact for our proof.
\begin{lem}\label{2.1}
Let $m\geq3$ and $u$ be a radial solution to \eqref{2.4}, if $a_{m-2}=\Delta^{m-2}u(0)=0$. Then $\lim_{r\rightarrow+\infty}\Delta^{m-1}u(r)<0$.
\end{lem}

\noindent
Proof. Let $v = \Delta^{m-1} u$. As $\Delta v =-e^u<0$,  $v(r)$ is decreasing in $\R_+$, so
$\lim_{r\rightarrow+\infty}v(r)= \ell \in \R\cup\{-\infty\}$ exists. Assume that $\ell\geq0$, then $v(r)>0$ in $\R_+$ and
$\Delta^{m-2}u(r)$ is increasing in $r$, which implies $\lim_{r\rightarrow+\infty}\Delta^{m-2}u(r)=\ell_1 > 0$, since
$\Delta^{m-2}u(0)=0$. By iterations, we conclude that
\begin{align*}
\lim_{r\rightarrow+\infty}\Delta^k u(r)= \infty,\quad \forall\; 0\leq k \leq m-3.
\end{align*}
Therefore $\lim_{r\to+\infty}\Delta v(r) = -\infty$. Again, by integrations, we get $\lim_{r\rightarrow+\infty}v(r)=-\infty$,
which contradicts $\ell \geq 0$, hence there holds $\ell < 0$.\qed

A useful consequence of Lemma \ref{2.1} is the following continuity result.
\begin{prop}\label{2.2}
Let $m \geq 3$ and $\Sigma_0 := \mathbb{R}^{m-2}\times{\{0}\}\times\mathbb{R}$. Then for any $(a_i)_{0\leq i\leq m-1}\in \Sigma_0$, i.e. $a_{m-2} = 0$, the radial solution $u_{(a_i)}$ to equation \eqref{2.4} satisfies
\begin{align*}
V(a_i) := \int_{\R^N} e^{u_{(a_i)}} dx <\infty.
\end{align*}
Moreover, the function $V$ is continuous in $\Sigma_0$.
\end{prop}

\noindent
Proof. Given $(a_i) \in \Sigma_0$, $\lim_{r\rightarrow+\infty}\Delta^{m-1}u_{(a_i)}(r)<0$ by Lemma \ref{2.1}, hence there is $R > 0$ large such that $\Delta^{m-1}u_{(a_i)}(R) < 0$.

\medskip
By ODE theory, the radial solution $u_{(a_i)}$ to equation \eqref{2.4} is continuous with respect to $(a_i)$ in $C^k_{loc}(\R^N)$ for any $k \in \N$. Consequently, there exists $\delta > 0$ small such that for any $|(a_i') -(a_i)| \leq \delta$, there holds
\begin{align*}
\left\|u_{(a_i')} - u_{(a_i)}\right\|_{C^{2m}(\overline B_{R})} \leq 1\quad \mbox{and} \quad \Delta^{m-1}u_{(a_i')}(R) \leq
\frac{\Delta^{m-1}u_{(a_i)}(R)}{2} := -M < 0.
\end{align*}
As $\Delta^{m-1}u$ is decreasing in $r$ for any radial solution to $\Delta^m u = -e^u$, we have $\Delta^{m-1}u_{(a_i')}(r) \leq -M$ if $r \geq R$ and $|(a_i') -(a_i)| \leq \delta$. Therefore, for $r \geq R$ and $|(a_i') -(a_i)| \leq \delta$,
\begin{align*}
& \Delta^{m-2}u_{(a_i')}(r) \\ = &\; \Delta^{m-2}u_{(a_i')}(R) +
\int_R^r\frac{1}{\rho^{N-1}}\left[R^{N-1}\left(\Delta^{m-2}u_{(a_i')}\right)'(R) + \int_R^\rho s^{N-1}\Delta^{m-1}u_{(a_i')}(s)ds\right] d\rho\\
\leq & \; C_1 + \int_R^r \left(-\frac{M}{N}\rho + C_2\rho^{1 - N}\right) d\rho\\
\leq & \; -\frac{Mr^2}{2N} + \frac{MR^2}{2N}+ C_1 + \frac{C_2R^{2-N}}{N-2}
= -\frac{Mr^2}{2N} + C_3.
\end{align*}
Here $C_i$ are some constants independent of $(a_i')$ verifying $|(a_i') -(a_i)| \leq \delta$. We get then $M' > 0$ and $R' \geq R$ such that
\begin{align*}
\Delta^{m-2}u_{(a_i')}(r) \leq -M' < 0, \quad \mbox{for all } r \geq R',\; |(a_i') -(a_i)| \leq \delta.
\end{align*}
By iterations, we can conclude that there exist $M_0 > 0$ and $R_0$ large such that
\begin{align*}
u_{(a_i')}(r) \leq -M_0r^{2m-4}, \quad \mbox{for all } r \geq R_0,\; |(a_i') -(a_i)| \leq \delta.
\end{align*}
Clearly $V(a_i) < \infty$ by the above estimate. It's not difficult to deduce the continuity of $V$ in $(a_i)$ using
the continuity of $u_{(a_i)}$ in $C^0_{loc}(\R^N)$ with respect to $(a_i)$, and the uniform estimate out of a compact set, we omit the details. \qed

\medskip
If $m=2$, we consider radial solutions to the following biharmonic equation
\begin{equation}\label{3.11}
\begin{cases}
\begin{aligned}
&\Delta^2 u=-e^u\\
&u'(0)=u'''(0)=0,\\
&\Delta u(0)= a, u(0)=-b.
\end{aligned}
\end{cases}
\end{equation}
Corresponding to Lemma \ref{2.1} for $m \geq 3$, we have
\begin{lem}\label{3.1}
For any $a, ~b\in\mathbb{R}$, the radial solution to \eqref{3.11} satisfies $\lim_{r\rightarrow+\infty} \Delta u(r)<0$.
\end{lem}

\noindent
Proof. Let $v = \Delta u$, as $\Delta v =-e^u<0$, $v$ is decreasing in $r \in(0,+\infty)$. So $\lim_{r\rightarrow+\infty} v(r)= \ell$
exists. If $\ell \geq0$, we have $v(r)>0$ in $\R_+$, then $u$ is increasing in $r$ and $\Delta v = -e^u \leq -e^{u(0)} = -e^{-b}$ in $\R^N$.
Then $\lim_{r\to+\infty} v(r) = -\infty$ since
$$
v(r) - a = \int_0^r\frac{1}{\omega_{N-1}\rho^{N-1}}\int_{B_\rho} \Delta vdx d\rho \leq -e^{-b}\frac{r^2}{2N}.
$$
This contradicts the assumption $\ell \geq 0$. So $\ell < 0$. \qed

Here and after, $\omega_{N-1}$ denotes the volume of the standard sphere $\S^{N-1} \subset \R^N$. Denote $u_{a, b}$ the radial solution to \eqref{3.11} and
$$\widetilde V(a,b):=\int_{\R^N} e^{u_{a, b}}dx.$$
Using Lemma \ref{3.1}, we can prove the continuity of $\widetilde V$ very similarly as for Proposition \ref{2.2}, so we omit the proof.
\begin{prop}\label{3.2}
For any $(a,b)\in\mathbb{R}^2$, then $\widetilde V(a,b)< \infty$. Moreover, $\widetilde V$ is continuous in $\R^2$.
\end{prop}

\subsection{Solutions with large volume for \eqref{1.1}} Here we prove the existence of radial solutions to \eqref{1.1}
with any large volume. As mentioned in Remark \ref{rem2.1}, the problem is relevant only in the conformal dimension. From now on, we fix $N = 2m$, even similar result holds true for any $N \geq 3$. The crucial point is to consider some special initial conditions.

\medskip
More precisely, for $m \geq 3$ and $N = 2m$, let
\begin{align}
\label{2.5}
c_0 := 4^{m-1} \times \prod^{m-1}_{k=1}\big[k(m-1+k)\big]
\end{align}
and consider \eqref{2.4} with $(a_i) = (-b, 0, \ldots 0, c_0) \in \Sigma_0 = \R^{m-2}\times\{0\}\times \R$. 
\begin{equation}\label{2.7}
\begin{cases}
\begin{aligned}
&\Delta^m u=-e^u, &&\\&u^{(2i+1)}(0)= \Delta^k u(0) = 0, &&  \forall\; i = 0\ldots m-1;\; k = 1, \ldots m-2,\\
&\Delta^{m-1} u(0) = c_0,&&\\
&u(0)=-b < 0.&&
\end{aligned}
\end{cases}
\end{equation}

\begin{thm}\label{2.3}
Let $m \geq 3$, denote $u_b$ the radial solution to equation \eqref{2.7}. Then
\begin{equation}\label{2.8}
\lim_{b\rightarrow+\infty}\int_{\R^{2m}} e^{u_b} dx =+\infty.
\end{equation}
Similarly, let $a = 8$ in \eqref{3.11} for $m=2$, there holds
\begin{equation}\label{3.12}
\lim_{b\rightarrow+\infty}\widetilde V(8, b) =+\infty.
\end{equation}
\end{thm}

\noindent Proof. We handle the cases $m \geq 3$ and $m = 2$ together. For simplicity and without confusion, we denote by $u$ the solution to \eqref{2.7} or
 the solution to \eqref{3.11} with $a = 8$.

\medskip
For any $m \geq 2$, let $\Phi(x) = \Phi(r)=r^{2m-2}-b$. Hence $\Delta^m \Phi=0$ in $\R^{2m}$ and $\Delta^i \Phi(0) = \Delta^iu(0)$
for any $0 \leq i \leq m-1$.  Set $w=u-\Phi,$ then $\Delta^m w= \Delta^m u =-e^u < 0$ and $\Delta^iw(0)=0$ for $0\leq i\leq m-1.$ By
iterations, we deduce easily that $\Delta^i w\leq 0$ in $\R^{2m}$ for $0\leq i \leq m-1$. In particular, $w \leq 0$ in $\R^{2m}$, i.e. $u\leq\Phi$ in $\R^{2m}$.

\medskip
Let $R_0 :=b^{\frac{1}{2m-2}}$, the unique zero of $\Phi$ in $(0, \infty)$. To prove \eqref{2.8} or \eqref{3.12}, we proceed by three steps.

\medskip
{\it Step 1}. Estimate of $\Delta^{m-1} w(R_0)$.

As $\Delta^m w=-e^u$, we have, for any $r > 0$,
$$
\begin{aligned}
 \Delta^{m-1}w(r) = - \int_0^r\frac{1}{\rho^{2m-1}}\int_0^\rho e^{u(s)}s^{2m-1} ds d\rho &\geq-\int_0^r \rho^{1-2m}\int_0^\rho e^{\Phi(s)} s^{2m-1}ds d\rho\\
 &=-\int_0^r e^{s^{2m-2}-b}s^{2m-1}ds\int_s^r \rho^{1-2m}d\rho\\
 &=-\frac{1}{2m-2}\int_0^r e^{s^{2m-2}-b}s\left[1-\left(\frac{s}{r}\right)^{2m-2}\right]ds\\
 &\geq-\int_0^r e^{s^{2m-2}-b}\left(1-\frac{s}{r}\right)sds\\
 & =-\frac{r^2}{2}\int_0^1 e^{r^{2m-2}t^{m-1}-b}\big(1-\sqrt{t}\big)dt.
\end{aligned}
$$
For the second inequality, we used the convexity of the function $h(x) = x^{2m-2}$ in $\R_+$ and we applied the change of variable $t = r^{-2}s^2$ for the last line. Therefore,
$$
\Delta^{m-1}w(R_0)\geq -\frac{R_0^2}{2}\int_0^1 e^{b(t^{m-1} - 1)}\big(1-\sqrt{t}\big)dt =:\xi(b).
$$
Moreover, there exists $\l > 0$ depending on $m$ such that $\l(1 - t) \leq 1 - t^{m-1}$ for any $t \in [0, 1]$. So we get
\begin{align*}
\int_0^1 e^{b(t^{m-1} - 1)}(1 - \sqrt{t})dt & \leq \int_0^1 e^{\l b(t - 1)}(1 - t)dt\\
& = \frac{1}{\l^2b^2} - \left(\frac{1}{\l b} +\frac{1}{\l^2b^2}\right) e^{-\l b} = O\left(\frac{1}{b^2}\right) \quad\text{as}~b\rightarrow\infty.
\end{align*}
As $\Delta^{m-1}w(r)$ is decreasing in $r$, there holds
\begin{align}
\label{new2}
\Delta^{m-1}w(r)\geq \Delta^{m-1}w(R_0)\geq \xi(b) = O\left(b^{\frac{1}{m-1} -2}\right) \quad \mbox{for } r \leq R_0.
\end{align}

{\it Step 2}. Estimates of $\Delta^i u(r)$, $i=0,...,m-1$ for $r> R_0$.

Define $r_0=\inf{\{r>0,~u(r)=0}\}\in (R_0, +\infty]$. We claim that
\begin{align}
\label{new1}
\lim_{b\rightarrow+\infty} (r_0-R_0)=0.
\end{align}
Remark that $u \leq 0$ in $[R_0, r_0]$, so
$\Delta^m w=-e^u\geq-1$ if $R_0 \leq r \leq r_0$. Therefore, if $r\in[R_0, r_0]$,
$$
\begin{aligned}
\Delta^{m-1}w(r)&\geq \Delta^{m-1}w(R_0)- \int_{R_0}^r\frac{1}{\omega_{2m-1}\rho^{2m-1}} \left(\int_{B_{R_0}} e^{\Phi} dx+\int_{B_\rho\setminus B_{R_0}} dx \right)d\rho,\\
&=\Delta^{m-1}w(R_0)- \int_{R_0}^r\rho^{1-2m}\left(\int_0^{R_0}e^{s^{2m-2}-b}s^{2m-1} ds +\int_{R_0}^\rho s^{2m-1}ds\right)d\rho.
\end{aligned}
$$
When $b\rightarrow+\infty$, there holds
\begin{align*}
 \eta(b) := \int_0^{R_0} e^{s^{2m-2}-b}s^{2m-1} ds = R_0^{2m} \int_0^1 e^{-b(1-t^{2m-2})}t^{2m-1} dt = O\left(\frac{R_0^{2m}}{b}\right).
\end{align*}
We obtain that for $r\in[R_0, r_0]$,
\begin{equation}\label{2.9}
\Delta^{m-1}w(r)\geq \Delta^{m-1}w(R_0)- \eta(b)R_0^{1-2m}(r-R_0)-\frac{1}{2m}\int_{R_0}^r\rho\left[1-\left(\frac{R_0}{\rho}\right)^{2m}\right] d\rho.
\end{equation}
On the other hand, by the convexity of $h(x) = x^{2m}$ in $\R_+$,
\begin{align*}
\frac{1}{2m}\int_{R_0}^r\rho\left[1-\left(\frac{R_0}{\rho}\right)^{2m}\right] d\rho \leq
\int_{R_0}^r \rho\left(1-\frac{R_0}{\rho}\right) d\rho = \frac{(r-R_0)^2}{2}, \quad \forall\; r \geq R_0.
\end{align*}

Denote $\widetilde{r_0} := \min{\{r_0, R_0 +1}\}$. Combining \eqref{new2} and \eqref{2.9}, for $r\in[R_0,\widetilde{r_0}]$, we have (as $m \geq 2$)
\begin{align}
\label{new3}
\Delta^{m-1}w(r) &\geq O\left(b^{\frac{1}{m-1}-2}\right) - O\left(\frac{R_0}{b}\right)(r - R_0) - \frac{(r-R_0)^2}{2} \geq O(1).
\end{align}
Using again \eqref{new2}, we obtain
$$
\Delta^{m-1}w(r)\geq O\left(b^{\frac{1}{m-1}-2}\right) +O(1)\chi_{[R_0, \widetilde{r_0}]}, \quad \forall \; r\in[0, \widetilde{r_0}].
$$
Here and in the following, $\chi_A$ denotes the characteristic function of a subset $A$ and $O(1)$ denotes a quantity uniformly bounded for $b$ sufficiently large.

\medskip
By iterations, for $0\leq j\leq m-2$ and $r\in[0, \widetilde{r_0}]$, we get
\begin{align}
\Delta^jw(r)\geq O\left(b^{\frac{1}{m-1}-2}\right)r^{2(m-1-j)} + O(1)\chi_{[R_0, \widetilde{r_0}]}.
\end{align}
In particular, let $j = 0$, there holds
$$
u(r)\geq \Phi(r)+ O\left(b^{\frac{1}{m-1}-2}\right)r^{2(m-1)} + O(1)\chi_{[R_0, \widetilde{r_0}]}, \quad \forall\; r\in[0, \widetilde{r_0}].
$$
Using the convexity of $\Phi$, we have then
$$
u(r)\geq(2m-2)R_0^{2m-3}(r-R_0)+ O\left(b^{\frac{1}{m-1}-2}\right)r^{2(m-1)} + O(1) \quad \mbox{in }\; [R_0, \widetilde{r_0}].
$$
Fix any $\varepsilon \in (0, 1)$, suppose that $r_0>R_0+\varepsilon$. Then
$$
0>u(R_0+\varepsilon)\geq(2m-2)b^{\frac{2m-3}{2m-2}}\varepsilon+ O\left(b^{\frac{1}{m-1}-1}\right)+ O(1),
$$
which is impossible for $b$ large enough, since
$$
\frac{2m-3}{2m-2} -\frac{1}{m-1} +1 = \frac{4m - 7}{2(m-1)} > 0, \quad \forall\; m \geq 2.
$$
In other words, when $b$ is sufficiently large, we have $r_0\leq R_0+\varepsilon$, so the claim \eqref{new1} is proved. An immediate consequence of \eqref{new1} is
\begin{equation}\label{2.10}
\liminf_{b\rightarrow+\infty}\Delta^{m-1}u(r_0)\geq c_0.
\end{equation}
Indeed, applying the first inequality in \eqref{new3},
$$
\Delta^{m-1}w(r_0)\geq O\left(b^{\frac{1}{m-1}-2}\right) - O\left(\frac{R_0}{b}\right)(r_0 - R_0) - \frac{(r_0-R_0)^2}{2},
$$
we get $\liminf_{b\rightarrow+\infty}\Delta^{m-1}w(r_0)\geq 0$ by \eqref{new1}, hence \eqref{2.10} holds true as $\Delta^{m-1}\Phi\equiv c_0$.

\medskip
{\it Step 3}. The proof of \eqref{2.8} and \eqref{3.12}.

Consider first $m \geq 3$. Recall that we denote by $u$, the radial solution of \eqref{2.7}. Let
$$V(b) :=\int_{\R^{2m}} e^u dx=\omega_{2m-1}\int_0^\infty e^{u(s)}s^{2m-1}ds.$$ By equation $\Delta^m u=-e^u$, we get
$$
r^{2m-1}(\Delta^{m-1} u)'(r)=-\int_0^r e^{u(s)}s^{2m-1} ds\geq-\frac{V(b)}{\omega_{2m-1}}.
$$
For any $r > r_0$, using the above inequality on $[r_0, r]$, there holds
\begin{align*}
\Delta^{m-1} u(r_0)\leq \Delta^{m-1} u(r)+\frac{V(b)}{\omega_{2m-1}}\int_{r_0}^r s^{1-2m}ds & \leq \Delta^{m-1} u(r)+\frac{V(b)}{\omega_{2m-1}}\int_{r_0}^\infty s^{1-2m}ds\\
& =\Delta^{m-1} u(r)+ \frac{V(b)r_0^{2-2m}}{(2m-2)\omega_{2m-1}}.
\end{align*}
Tending $r$ to $+\infty$, we conclude by Lemma \ref{2.1}
that $$
V(b)>(2m-2)\omega_{2m-1}r_0^{2m-2}\Delta^{m-1} u(r_0).
$$
Hence $\lim_{b\rightarrow+\infty} V(b)=+\infty$ by \eqref{2.10} and $\lim_{b\to +\infty} r_0 = +\infty$.

\medskip
The proof of \eqref{3.12} is completely similar, so we omit it.\qed

\begin{rem}
The formula \eqref{3.12} gives a positive answer to a question in \cite{M}, page 981. Assume that $u$ solves \eqref{3.11} with $a = 8$ in $\R^4$. Let $v(x) = u(\l x) + 4\ln\l$ with $\l = e^{b/4}$, then $\Delta^2 v = -e^v$, $v(0) = 0$ and $v''(0) = e^{b/2}u''(0) = 2e^{b/2}$. Hence $v''(0) \to +\infty$ is equivalent to $b \to +\infty$.
\end{rem}

\subsection{Proof of Theorem \ref{thm1.1} completed} Let $m \geq 3$ and $\widetilde u_b$ be the radial solution to equation \eqref{2.4} with $(a_i) = (-b,0,...0,-c_0) \in \Sigma_0$. As above, there holds $\widetilde u_b \leq \Psi(r) := -r^{2m-2} -b$ in $\R^{2m}$. Hence
$$
\lim_{b \rightarrow+\infty}\int_{\R^{2m}} e^{\widetilde u_b} dx \leq\lim_{b \rightarrow+\infty}\int_{\R^{2m}} e^{\Psi} dx = \lim_{b \rightarrow+\infty}\int_{\R^{2m}} e^{-|x|^{2m-2}-b}dx=0.
$$
By Proposition \ref{2.2}, Theorem \ref{2.3} and the above estimate, we get readily that $V(\Sigma_0) = (0, \infty)$, so we are done.

\medskip
The argument for $m = 2$ is completely similar. Considering the radial solution $\widetilde u_b$ to \eqref{3.11} with $a = -8$, we prove easily that $\inf_{\R^2} \widetilde V(a, b) = 0$. Using \eqref{3.12} and Proposition \ref{3.2}, there holds $\widetilde V(\R^2) = (0, \infty)$.\qed

\section{Proof of Theorem \ref{thm1.1new}}
\reset
For \eqref{1.1new}, we use a different approach, which is based on the following well-known comparison result (see for instance Proposition 13.2 in \cite{FF})
\begin{lem}
\label{comp}
Let $u, v\in C^{2m}([0, R))$ be two radial functions such that $\Delta^m u - e^u \geq \Delta^m v - e^v$ in $[0, R)$ and
\begin{align}
\Delta^ku(0) \geq \Delta^kv(0), \; (\Delta^ku)'(0) \geq (\Delta^kv)'(0), \quad \forall \; 0\leq k \leq m-1.
\end{align}
Then we have $u \geq v$ in $[0, R)$.
\end{lem}

Let $u_0$ be an entire radial solution of \eqref{1.1new} with $V_0 := \|e^{u_0}\|_{L^1(\R^{2m})} \in (0, \infty]$, consider $u_\alpha$
the solution to the following initial value problem
\begin{equation}\label{4.1}
\begin{cases}
\begin{aligned}
&\Delta^m u= e^u, && \\
&u^{(2i+1)}(0)=0, && \forall\; 0\leq i\leq m-1,\\
&\Delta^iu(0)= \Delta^iu_0(0), && \forall\; 0\leq i\leq m-2,\\
& \Delta^{m-1} u(0) = \Delta^{m-1} u_0(0) - \alpha, && \alpha > 0.
\end{aligned}
\end{cases}
\end{equation}
For any $\alpha > 0$, by Lemma \ref{comp}, we have $u_\alpha \leq u_0$ whenever it exists. On the other hand, we have $u_\alpha \geq
\Phi_\alpha = \sum_{0\leq j\leq m-1} \alpha_jr^{2j}$ with $\alpha_j \in \R$ verifying $\Delta^i\Phi_\alpha(0) = \Delta^iu_\alpha(0)$ for
$i = 0, \ldots, m-1$. Then no blow-up will occur for $u_\alpha$ in any compact set, which means that
$u_\alpha$ is globally defined and $\alpha \mapsto u_\alpha$ is a decreasing family of functions in $\R^{2m}$ by Lemma \ref{comp}. We claim that
\begin{align}
\label{4.2}
\lim_{\alpha \to \infty} u_\alpha(r) = -\infty, \quad \mbox{for any } r \geq 0.
\end{align}
Let $v_\alpha = \Delta^{m-1}u_\alpha$, as $0 \leq \Delta v_\alpha = e^{u_\alpha} \leq e^{u_0}$ and
\begin{align*}
v_\alpha(r) = \Delta^{m-1}u_0(0) - \alpha + \frac{1}{2m-2}\int_0^r \Delta v_\alpha(s)s\left[1 - \left(\frac{s}{r}\right)^{2m-2}\right] ds,
\end{align*}
we get readily that $v_\alpha$ tends uniformly to $-\infty$ in any compact set of $\R_+$ as $\alpha \to \infty$. By iterations,
we obtain that $\Delta^{m-2}u_\alpha, \ldots
\Delta u_\alpha$ tend to $-\infty$ uniformly in any compact set of $\R_+$, hence \eqref{4.2} is satisfied.

\medskip
Moreover, by Lemmas 7.6 and 7.8 in \cite{FF}, for any $\alpha > 0$, there holds $\lim_{r\to \infty}\Delta^{m-1} u_\alpha(r) < 0$ hence
$u_\alpha(r) \leq -C_\alpha r^{2m-2}$ for $r$ large enough with some $C_\alpha > 0$. Therefore $e^{u_\alpha} \in L^1(\R^{2m})$ for
any $\alpha > 0$, and $\alpha
\mapsto \|e^{u_\alpha} \|_{L^1(\R^{2m})}$ is continuous in $(0, \infty)$, combining with the monotonicity w.r.t.~$\alpha$. Furthermore, the claim
\eqref{4.2} implies then
\begin{align*}
 \lim_{\alpha \to \infty}\int_{\R^{2m}} e^{u_\alpha} dx = 0.
\end{align*}
As the monotonicity of $u_\alpha$ gives also
\begin{align*}
 \lim_{\alpha \to 0^+}\int_{\R^{2m}} e^{u_\alpha} dx = \int_{\R^{2m}} e^{u_0} dx = V_0,
\end{align*}
the proof is completed. \qed

\section{Further remarks and open questions}
By the proof of Proposition \ref{2.2}, for $m \geq 3$ and $(a_i) \in \R^m$, if the solution of \eqref{2.4} verifies $\lim_{r\to+\infty}\Delta^{m-1}u_{(a_i)}(r) < 0$, then vol$(g_u)$ is finite and the function $V$ with the initial data as variables is
continuous at the point $(a_i)$. However we have no answer for the following question.

\medskip
\noindent
{\bf Question 1}: For any $(a_i) \in \R^m$ with $m \geq 3$, let $u$ be the solution of \eqref{2.4}, is the total volume vol$(g_u)$ finite?
If the answer is yes, is the volume function $V$ continuous in whole $\R^{2m}$?

\medskip
Another natural question comes from Theorem \ref{even}. Our approach is to study radial solutions of \eqref{1.1new}. By Hyder \& Martinazzi's result in \cite{HM} on the negative constant $Q$-curvature situation, for $m \geq 3$ {\it odd}, the radial solutions of \eqref{1.1new} can provide arbitrary volume. However, this is not always true for $m$ {\it even}, since when $m = 2$, ${\mathcal V} = (0, {\rm vol}(\S^4)]$ by \cite{L, CC, WY}.

\medskip
\noindent
{\bf Question 2}: Let $m$ be even and $m \geq 4$, do we have ${\mathcal V} = (0, \infty)$ for radial solutions of \eqref{1.1new}?

\medskip
Consider the radial solutions to \eqref{1.1new} as a initial value problem with $\Delta^k u(0) = \beta_k$, $0 \leq k \leq m-1$. By Theorem 2.2 in \cite{FF}, there exists a function $\Phi: \R^{m-2} \rightarrow (-\infty, 0)$ such that the solution $u$ is globally defined in $\R^{2m}$ if and only if $\beta_{m-1} \leq \Phi(\beta_1,\ldots \beta_{m-2})$. On the other hand, for given $\beta_1,\ldots \beta_{m-2}$, the solution $u_{(\beta_i)}$ is increasing w.r.t.~$\beta_{m-1}$ by Lemma \ref{comp}, so \begin{align*}
\sup{\mathcal V} = \sup\left\{\int_{\R^{2m}} e^{2mu}dx, \;\; u = u_{(\beta_i)} \mbox{ with } \beta_{m-1} = \Phi(\beta_1,\ldots \beta_{m-2}).\right\}.
\end{align*}
Therefore, to answer the above question, we need just to understand the radial solutions with $(\beta_i)$ on the boundary hypersurface for the global existence. Unfortunately we have few information for these borderline entire radial solutions. For instance, we don't know the asymptotic decay of such solutions as $r \to \infty$, see Theorem 2.5 and Problem 2.1 (ii) in \cite{FF}.

\medskip
A last question concerns the infinite volume entire solutions. When $m = 1$, Liouville proved that given a holomorphic function $h$ in $\O \subset \C$, the function
\begin{align*}
u(z) := \ln\frac{2|h'(z)|}{1 + |h(z)|^2}
\end{align*}
satisfies $-\Delta u = e^{2u}$ in $\O\backslash\{z\in \O, h'(z) = 0\}$. The conformal metrics in $\R^2$ with vol$(g_u) < \infty$, i.e.~the solutions to \eqref{02} correspond to $h(z) = az + b$ with $a, b \in \C$. So we can describe many entire solutions of $-\Delta u = e^{2u}$ in $\R^2$ with infinite volume.

\medskip
For $m \geq 2$, of course we can use entire radial solutions $v$ of $(-\Delta)^m v = e^{2mv}$ in $\R^N$ with $3 \leq N \leq 2m - 1$ to construct constant $Q$-curvature conformal metrics in $\R^{2m}$ with infinite volume, for example by considering $u(x) := v(x_1,\ldots x_N)$. However, we wonder if other examples exist.

\medskip
{\bf Question 3}: For $m \geq 2$, are there entire solutions of $(-\Delta)^m u = e^{2mu}$ in $\R^{2m}$ such that $u$ does not allow any symmetry and $e^{2mu} \not\in L^1(\R^{2m})$?

\medskip
\noindent
{\bf Acknowledgements.}
This work is realized during the visit of Huang at the Institut Elie Cartan de Lorraine. She would like to thank the institute for 
its warm hospitality, and the China Scholarship Council for supporting this visit in Metz. Huang is also partially supported by 
NSFC (No.~11271133).

\end{document}